\documentclass{article}
\usepackage{amsmath,amsthm,amssymb}

{\theoremstyle{definition}\newtheorem{definition}{Definition}[section]
\newtheorem{notation}[definition]{Notation}

\newtheorem{remark}[definition]{Remark}

}

\newtheorem{proposition}[definition]{Proposition}
\newtheorem{lemma}[definition]{Lemma}
\newtheorem{theorem}[definition]{Theorem}

\newtheorem{corollary}[definition]{Corollary}

\newcommand{\Ktil}{\widetilde{K}}
\newcommand{\Hincl}{H_{\text{\rm incl}}}
\newcommand{\Hrep}{H_{\text{\rm rep}}}
\newcommand{\FAlg}{\operatorname{FAlg}}

\newcommand{\si}{\sigma}
\newcommand{\recht}{\rightarrow}
\newcommand{\Ker}{\operatorname{Ker}}
\newcommand{\Out}{\operatorname{Out}}
\newcommand{\Aut}{\operatorname{Aut}}

\newcommand{\cZ}{\mathcal{Z}}
\newcommand{\actson}{\curvearrowright}

\newcommand{\al}{\alpha}

\newcommand{\Ad}{\operatorname{Ad}}
\newcommand{\om}{\omega}

\newcommand{\Z}{\mathbb{Z}}
\newcommand{\cU}{\mathcal{U}}
\newcommand{\SL}{\operatorname{SL}}

\newcommand{\be}{\beta}

\newcommand{\id}{\operatorname{id}}

\newcommand{\C}{\mathbb{C}}
\newcommand{\cF}{\mathcal{F}}
\newcommand{\eps}{\varepsilon}
\newcommand{\cC}{\mathcal{C}}

\newcommand{\cL}{\mathcal{L}}

\newcommand{\M}{\operatorname{M}}
\newcommand{\ot}{\otimes}
\newcommand{\N}{\mathbb{N}}

\newcommand{\notembed}[1]{\underset{#1}{\not\prec}}
\newcommand{\embed}[1]{\underset{#1}{\prec}}
\newcommand{\fembed}[1]{\underset{#1}{\overset{f}{\prec}}}

\newcommand{\Tr}{\operatorname{Tr}}
\newcommand{\dis}{\displaystyle}

\newcommand{\R}{\mathbb{R}}

\newcommand{\QN}{\operatorname{QN}}
\newcommand{\dpr}{^{\prime\prime}}
\newcommand{\la}{\langle}
\newcommand{\ra}{\rangle}
\newcommand{\cN}{\mathcal{N}}
\newcommand{\cV}{\mathcal{V}}

\newcommand{\cH}{\mathcal{H}}
\newcommand{\otalg}{\ot_{\text{\rm alg}}}
\newcommand{\proj}{p^{K(\al)}_{\text{\rm centr}}}
\newcommand{\cW}{\mathcal{W}}

\begin{document}
\begin{center}
{\LARGE\bf Factors of type II$_1$ without non-trivial \vspace{2mm} \\ finite index subfactors}

\bigskip

{\sc by Stefaan Vaes}

{\small Department of Mathematics; K.U.Leuven; Celestijnenlaan 200B; B--3001 Leuven (Belgium) \\ E-mail: stefaan.vaes@wis.kuleuven.be}
\end{center}

\begin{abstract}
\noindent We call a subfactor $N \subset M$ trivial if it is isomorphic with the obvious inclusion of $N$ in $\M_n(\C) \ot N$. We prove the existence of type II$_1$ factors $M$ without non-trivial finite index subfactors. Equivalently, every $M$-$M$-bimodule with finite coupling constant, both as a left and as a right $M$-module, is a multiple of $L^2(M)$. Our results rely on the recent work of
Ioana, Peterson and Popa, who proved the existence of type II$_1$ factors without outer automorphisms.
\end{abstract}

\section*{Introduction}

We say that a subfactor $N \subset M$ of finite index is \emph{trivial}, if there exists $n \in \N$ such that $N \subset M$ is isomorphic with $1 \ot N \subset \M_n(\C) \ot N$. We prove that there exist type II$_1$ factors all of whose finite index subfactors are trivial. An $M$-$M$-bimodule $_M H_M$ is said to be \emph{bifinite} if $\dim(H_M) < \infty$ and $\dim(_M H) < \infty$.
In the language of Connes' correspondences, our main theorem then tells that there exist type II$_1$ factors $M$ such that every bifinite $M$-$M$-bimodule is trivial, i.e.\ isomorphic with a direct sum of copies of $_M L^2(M)_M$.

Such II$_1$ factors are very special. Indeed, any automorphism $\al \in \Aut(M)$ gives rise to an $M$-$M$-bimodule $H(\al)$ on the Hilbert space $L^2(M)$ by the formula
$$x \cdot \xi = \al(x) \xi \quad\text{and}\quad \xi \cdot x = \xi x \quad\text{for all}\;\; x \in M, \xi \in L^2(M) \; .$$
This $M$-$M$-bimodule is trivial if and only if $\al$ is an inner automorphism. So, absence of non-trivial finite index subfactors implies \emph{absence of outer automorphisms}. Further, if $p$ is a projection in $M$ and $\pi : M \recht pMp$ a $^*$-isomorphism, one considers analogously the $M$-$M$-bimodule $_{\pi(M)} L^2(pM)_M$. Hence, absence of non-trivial finite index subfactors implies \emph{triviality of the fundamental group}.

Because of the constructions in the previous paragraph, the bifinite $M$-$M$-bimodules, should be considered as the \emph{generalized symmetries} of the II$_1$ factor $M$. Our main statement then becomes that there exist type II$_1$ factors all of whose generalized symmetries are inner.

In general, computing the outer automorphism group $\Out(M)$ of a II$_1$ factor $M$ is very hard. Connes discovered in \cite{C1} that $\Out(M)$ is countable whenever $M$ is the group von Neumann algebra of an ICC property (T) group. Only very recently, Ioana, Peterson and Popa proved the existence of type II$_1$ factors $M$ with $\Out(M)$ trivial, see \cite{IPP}. Their theorem is an existence result in the same way as is the main result in the current paper. We comment on that below. Explicit examples of II$_1$ factors with trivial outer automorphism group were constructed by Popa and the author in \cite{PV}, using crossed products by generalized Bernoulli actions and relying on the techniques of Popa's breakthrough von Neumann strong rigidity results in \cite{P1,P2}. Note that in \cite{PV}, it is shown as well that any group of finite presentation can be explicitly realized as the outer automorphism group of a II$_1$ factor.

Also the fundamental group of a II$_1$ factor, introduced by Murray and von Neumann in \cite{MvN}, is very hard to compute, unless, of course, you deal with a McDuff factor and get $\R^*_+$ as its fundamental group. Connes proved in \cite{C1} that the fundamental group of the group von Neumann algebra of an ICC property (T) group is countable. The first example of a II$_1$ factor with trivial fundamental group was given by Popa in \cite{P5}, as the group von Neumann algebra of $\SL(2,\Z) \ltimes \Z^2$. Many other such examples are given in \cite{IPP,P1,P2,PV}. In \cite{P1}, Popa constructs type II$_1$ factors with an arbitrarily prescribed countable subgroup of $\R^*_+$ as a fundamental group. An alternative construction is given in \cite{IPP}.

The type II$_1$ factors studied in this paper are of the form $M=R \rtimes \Gamma$, where $\Gamma = \Gamma_0 * \Gamma_1$ is the free product of two infinite groups and $\Gamma \actson R$ is an action by outer automorphisms on the hyperfinite II$_1$ factor $R$. We formulate strong conditions on the groups and the actions involved, that ensure that all bifinite $M$-$M$-bimodules, are trivial. But, we do not give explicit examples of actions that satisfy all these requirements: as in \cite{IPP}, we rather prove the existence of such actions through a Baire category argument.

The following argument, due to Ioana, Peterson and Popa \cite{IPP} is a key ingredient to prove that, under suitable conditions, every bifinite $M$-$M$-bimodule is trivial when $M = R \rtimes (\Gamma_0 * \Gamma_1)$. One first assumes that $R \subset M$ has the relative property (T). The free product $\Gamma_0 * \Gamma_1$ gives rise to a strong deformation property of $M$. Combined with the relative property (T) for $R \subset M$, this fixes somehow the position of $R$ inside $M$. More precisely, it allows to conclude that any finite index inclusion $\pi : M \recht M^t$ can be unitarily conjugated into one in which $\pi(R) \subset R^t$, see Theorem \ref{thm.crucialstep} and Propositions \ref{prop.propone} and \ref{prop.proptwo}.

\subsection*{Acknowledgment}

The author is most grateful to Dietmar Bisch, Adrian Ioana, Jesse Peterson and Sorin Popa for the many fruitful discussions and their kind hospitality at Vanderbilt University and at the University of California at Los Angeles.

\section*{Notations, terminology and preliminaries}

Throughout, $(M,\tau)$ denotes a von Neumann algebra $M$ with a faithful normal tracial state $\tau$. We denote, for all $n \in \N_0$ and all $(M,\tau)$,
$$M^n := \M_n(\C) \ot M \; .$$
We use the convention $\N_0 = \{1,2,\ldots\}$.
If $M$ is a II$_1$ factor and $t > 0$, we also introduce the usual notation $M^t = p M^n p$ whenever $p \in M^n$ is a projection with non-normalized trace equal to $t$.

We make an extensive use of Popa's technique of intertwining subalgebras using bimodules. We explain a few notations and refer to the Appendix for more detailed statements. Let $(M,\tau)$ be a von Neumann algebra with a fixed faithful normal tracial state $\tau$. Let $A, B \subset M$ be von Neumann subalgebras. We say that $A$ embeds into $B$ inside $M$ and write
$$A \embed{M} B$$ if $L^2(M)$ contains a non-zero $A$-$B$-subbimodule $H$ that is finitely generated as a right $B$-module. We write $$A \fembed{M} B$$ if for every non-zero projection $p \in A' \cap M$, $L^2(pM)$ contains a non-zero $A$-$B$-subbimodule that is finitely generated as a right $B$-module.

The \emph{normalizer} of $A \subset M$ consists of the unitaries $u \in \cU(M)$ satisfying $uAu^* = A$ and is denoted by $\cN_M(A)$. We say that $A \subset M$ is \emph{regular} if $\cN_M(A)^{\prime\prime} = M$.

If $A \subset (M,\tau)$ is a von Neumann subalgebra, we say that $a \in M$ \emph{quasi-normalizes} $A$ if there exist $a_1,\ldots,a_n, b_1,\ldots,b_m \in M$ satisfying $A a \subset \sum_{i=1}^n a_i A$ and $a A \subset \sum_{j=1}^m A b_j$. The set of elements quasi-normalizing $A$ is denoted by $\QN_M(A)$ and is a unital $^*$-subalgebra of $M$ containing $A$. We call \emph{quasi-normalizer} of $A$ inside $M$ the von Neumann algebra $\QN_M(A)\dpr$ generated by the elements quasi-normalizing $A$. If $\QN_M(A)\dpr = M$, we say that the inclusion $A \subset M$ is \emph{quasi-regular}.

If $A \subset (M,\tau)$ is a von Neumann subalgebra, Jones' \emph{basic construction} \cite{Jon3} is denoted by $\la M,e_A \ra$ and defined as the von Neumann algebra acting on $L^2(M)$ generated by $A$ and the orthogonal projection $e_A$ of $L^2(M)$ onto $L^2(A)$. Note that $A$ commutes with $e_A$ and that $e_A x e_A = E_A(x) e_A$ for all $x \in M$, where $E_A : M \recht A$ denotes the unique $\tau$-preserving conditional expectation. Equivalently, $\la M,e_A \ra$ equals the commutant of the right $A$-action on $L^2(M)$.

If $(A,\tau)$ is a von Neumann algebra with a fixed faithful normal tracial state $\tau$ and if $H_A$ is a right $A$-module, the commutant $A'$ of the right $A$-action on $H$ is equipped with a canonical normal faithful \emph{semifinite trace} $\Tr$ that can be characterized as follows:
$$\Tr(T T^*) = \tau(T^*T) \quad\text{whenever}\quad T : L^2(A) \recht A : T(\xi a) = (T \xi) a \;\;\text{for all}\;\; \xi \in H, a \in A \; .$$
One defines
$$\dim(H_A) := \Tr(1)$$
and one calls $\dim(H_A)$ the \emph{coupling constant} or the \emph{relative dimension} of the right $A$-module $H_A$. As such, the definition of $\dim(H_A)$ depends on the choice of tracial state $\tau$ on $A$. Throughout this paper, either $A$ will be a II$_1$ factor, in which case the coupling constant is canonically defined, or $A$ will inherit a trace from a natural ambient II$_1$ factor.

For II$_1$ factors, the coupling constant is canonically defined and it is then a \emph{complete invariant} of Hilbert $A$-modules. If $A$ has a non-trivial center, a complete invariant of Hilbert $A$-modules can be given in terms of the center-valued trace. We shall only use the following corollary: if $\dim(H_A) < \infty$ and $\eps > 0$, there exists a central projection $z \in \cZ(A)$, $n \in \N$ and a projection $p \in A^n$ such that $\tau(1-z) < \eps$ and $(Hz)_A \cong (p L^2(A)^{\oplus n})_A$ as $A$-modules.

Let $A \subset (M,\tau)$. Regarding the basic construction $\la M,e_A \ra$ as the commutant of the right $A$-action on $L^2(M)$, we get a natural normal faithful semifinite trace $\Tr$ on $\la M,e_A \ra$. It is characterized by the formula $\Tr(x e_A y) = \tau(xy)$, for all $x,y \in M$.

If $_M H_M$ is an $M$-$M$-bimodule and $A \subset M$ a von Neumann subalgebra, a vector $\xi \in H$ is said to be \emph{$A$-central} if $a \xi = \xi a$ for all $a \in A$.

In \cite{P5}, Popa defined the \emph{relative property (T)} for an inclusion $A \subset (M,\tau)$ of a von Neumann algebra $A$ into the von Neumann algebra $M$ equipped with a faithful normal tracial state $\tau$. An equivalent form of this definition goes as follows. For every $\eps > 0$, there exists a finite subset $\cF \subset M$ and a $\delta > 0$ such that every $M$-$M$-bimodule that admits a unit vector $\xi$ with the property
$$|\la \xi, a \xi b \ra - \tau(ab) | < \delta \quad\text{for all}\;\; a,b \in \cF \; ,$$
admits an $A$-central vector $\xi_0$ satisfying $\|\xi_0 - \xi \| < \eps$.

\section*{The fusion algebra of a II$_1$ factor}

If $M$ is a type II$_1$ factor and $_M H_M$ an $M$-$M$-bimodule, we say that $H$ is \emph{bifinite} if $\dim( _M H) < \infty$ and $\dim(H_M) < \infty$. The \emph{fusion algebra} of $M$ is defined as the set of all bifinite $M$-$M$-bimodules modulo isomorphism of bimodules and is denoted as $\FAlg(M)$. Note that $\FAlg(M)$ is equipped with the operations of direct sum and \emph{Connes tensor product}, see V.Appendix~B in \cite{C4} and the brief review below. One has the obvious notion of an \emph{irreducible element} in $\FAlg(M)$, and every element in $\FAlg(M)$ is the direct sum of a finite number of irreducibles.

Every $M$-$M$-bimodule $_M H_M$ has a \emph{contragredient} $M$-$M$-bimodule $_M \overline{H}_M$. Its carrier Hilbert space is the adjoint Hilbert space $\overline{H}$ while its bimodule structure is given by
$$x \cdot \overline{\xi} = \overline{\xi a^*} \quad\text{and}\quad \overline{\xi} \cdot a = \overline{a^* \xi} \; .$$
If $H$ and $K$ are bifinite $M$-$M$-bimodules, then $H$ and $K$ are disjoint if and only if $\overline{H} \ot_M K$ is disjoint from the trivial bimodule $_M L^2(M)_M$ if and only if $H \ot_M \overline{K}$ is disjoint from the trivial bimodule.

Finally, recall \emph{Frobenius reciprocity}: if $H, K, L \in \FAlg(M)$, the multiplicity of $H$ in $K \ot_M L$ equals the multiplicity of $K$ in $H \ot_M \overline{L}$ and equals the multiplicity of $L$ in $\overline{K} \ot_M H$.

We briefly recall the Connes tensor product. If $_M H_M$ is an $M$-$M$-bimodule, there is a natural dense subbimodule $\cH \subset H$ and $\cH$ is a W$^*$-$M$-$M$-bimodule, meaning that there is an $M$-valued scalar product on $\cH$. More precisely, $\cH$ consists of those vectors $\xi \in H$ such that there exists $\lambda \geq 0$ satisfying $\|\xi a\| \leq \lambda \|a\|_2$ for all $a \in M$. If now $_M K_M$ is another $M$-$M$-bimodule, the Connes tensor product $H \ot_M K$ is defined as the separation and completion of the algebraic tensor product $\cH \otalg K$ for the scalar product
$$\la a \ot \xi, b \ot \eta \ra := \la \xi, \; \la a, b \ra_M \; \eta \ra \; .$$
The $M$-$M$-bimodule structure on $H \otalg K$ is given by
$$a \cdot (b \ot \xi) = ab \ot \xi \quad\text{and}\quad (b \ot \xi) \cdot a = b \ot (\xi a) \; .$$

When there is no risk for misunderstanding, the tensor product $H \ot_M K$ is sometimes simply denoted by $HK$.

\begin{notation} \label{not.Hpsi}
Whenever $\psi : M \recht p M^n p$ is a finite index inclusion, we denote by $H(\psi)$ the bifinite $M$-$M$-bimodule with carrier Hilbert space $p(\M_{n,1}(\C) \ot L^2(M))$ and bimodule structure $x \cdot \xi = \psi(x)\xi$, $\xi \cdot x = \xi x$. In particular, every automorphism $\al \in \Aut(M)$ defines the element $H(\al) \in \FAlg(M)$ and as such, one considers $\Out(M) \subset \FAlg(M)$.
\end{notation}

Note that every bifinite $M$-$M$-bimodule is isomorphic with some $H(\psi)$. Moreover, if $\psi : M \recht p M^n p$ and $\theta : M \recht q M^m q$ are finite index inclusions, the $M$-$M$-bimodules $H(\psi)$ and $H(\theta)$ are isomorphic if and only if there exists a unitary $u \in p (\M_{n,m}(\C) \ot M) q$ satisfying $\theta(x) = u^* \psi(x) u$ for all $x \in M$. Also note that $H(\psi) \ot_M H(\theta) \cong \, H((\id \ot \theta)\psi)$.

A subset $\cF \subset \FAlg(M)$ is called a \emph{fusion subalgebra} if $\cF$ is closed under taking submodules, direct sums and tensor products. An important role is played in this paper by \emph{freeness} between fusion subalgebras.

\begin{definition} \label{def.free}
Let $M$ be a II$_1$ factor. Two fusion subalgebras $\cF_1,\cF_2 \subset \FAlg(M)$ are said to be \emph{free} if the following two conditions hold.
\begin{itemize}
\item Every tensor product of non-trivial irreducible bimodules, with factors alternatingly from $\cF_1$ and $\cF_2$, is irreducible.
\item Two tensor products of non-trivial irreducible bimodules, with factors alternatingly from $\cF_1$ and $\cF_2$, are equivalent if and only if they are factor by factor equivalent.
\end{itemize}
Equivalently, $\cF_1$ and $\cF_2$ are free if every tensor product of non-trivial irreducible bimodules, with factors alternatingly from $\cF_1$ and $\cF_2$, is disjoint from the trivial bimodule.
\end{definition}

Whenever $\al \in \Aut(M)$, we defined in Notation \ref{not.Hpsi} the bimodule $H(\al) \in \FAlg(M)$. So, if $\Gamma \actson M$ is an outer action, we can regard $\Gamma$ as a fusion subalgebra of $\FAlg(M)$.

\begin{definition}\label{def.almostnorm}
Let the countable group $\Gamma$ act outerly on the II$_1$ factor $N$. The \emph{almost normalizer of $\Gamma \actson N$ inside $\FAlg(N)$} is defined as the fusion subalgebra of $\FAlg(N)$ generated by the bifinite $N$-$N$-bimodules that can be realized as an $N$-$N$-subbimodules of a bifinite $(N \rtimes \Gamma)$-$(N \rtimes \Gamma)$-bimodule.
\end{definition}

We prove some results on the almost normalizing bimodules for $\Gamma \actson N$ in Section~\ref{sec.bifinite}. There, the terminology of \emph{bimodules almost normalizing} $\Gamma \actson N$, will become more clear as well. Right now, we already make the following observation.

\begin{lemma}\label{lem.almostnorm}
Let $\Gamma \actson N$ be an outer action on the II$_1$ factor $N$. If $\Gamma_0 < \Gamma$ is a finite index subgroup, the almost normalizers of $\Gamma_0 \actson N$ and $\Gamma \actson N$ inside $\FAlg(N)$, coincide.
\end{lemma}
\begin{proof}
Tensoring with the obvious inclusion bimodule $$\Hincl(\Gamma_0,\Gamma) = \, _{N \rtimes \Gamma_0} L^2(N \rtimes \Gamma)_{N \rtimes \Gamma}$$ and its contragredient, one goes back and forth between bifinite bimodules for $N \rtimes \Gamma_0$ and $N \rtimes \Gamma$.
\end{proof}

\section{Statement of the main result}

\begin{theorem}\label{thm.exists}
There exist II$_1$ factors $M$ with trivial fusion algebra: every bifinite $M$-$M$-bimodule is isomorphic with $_M (L^2(M)^{\oplus n})_M$ for some $n \in \N$.

In particular, $M$ has no outer automorphisms, has trivial fundamental group and only has trivial finite index subfactors: if $N \subset M$ is a finite index subfactor, $(N \subset M) \cong (1 \ot N \subset \M_n(\C) \ot N)$ for some $n \in \N$. In particular, every finite index irreducible subfactor of $M$, equals $M$.
\end{theorem}

The II$_1$ factors in the above theorem are of the form $M = R \rtimes \Gamma$, where $R$ is the hyperfinite II$_1$ factor, $\Gamma$ is the free product of two groups without non-trivial finite dimensional unitary representations and the outer action $\Gamma \actson N$ satisfies the following specific conditions.

\begin{theorem}\label{thm.main}
Let $\Gamma_0$, $\Gamma_1$ be infinite groups acting outerly on the II$_1$ factor $N$. Make the following assumptions.
\begin{itemize}
\item The groups $\Gamma_0, \Gamma_1, \Z$ are two by two not virtually isomorphic.
\item The groups $\Gamma_0, \Gamma_1$ are not virtually isomorphic to a non-trivial free product.
\item Denote by $\cF$ the fusion subalgebra of $\FAlg(N)$ consisting of the bifinite $N$-$N$-bimodules that almost normalize $\Gamma_0 \actson N$. Then, $\cF$ and $\Gamma_1$ are free as fusion subalgebras of $\FAlg(N)$. (See Defs.\ \ref{def.free} and \ref{def.almostnorm} for relevant terminology.)
\item $N \subset N \rtimes \Gamma_0$ has the relative property (T).
\end{itemize}
Set $M = N \rtimes (\Gamma_0 * \Gamma_1)$. If $_M H_M$ is a bifinite $M$-$M$-bimodule, there exists a finite-dimensional unitary representation $\theta : \Gamma_0 * \Gamma_1 \recht \cU(n)$, such that $_M H_M$ is isomorphic with the $M$-$M$-bimodule $\Hrep(\theta)$ defined below.
\end{theorem}

The $M$-$M$-bimodule $\Hrep(\theta)$ is defined as follows. The Hilbert space is given by $\C^n \ot L^2(M)$ and
$$(x u_g) \cdot \xi = (\theta(g) \ot x u_g) \xi \quad\text{and}\quad \xi \cdot y = \xi(1 \ot y)$$
for all $\xi \in \C^n \ot L^2(M)$, $g \in \Gamma_0 * \Gamma_1$, $x \in N$ and $y \in M$.

{\bf Organization of the proof of Theorems \ref{thm.exists} and \ref{thm.main}.} A given bifinite $M$-$M$-bimodule is of the form $H(\psi)$, where $\psi : N \rtimes \Gamma \recht (N \rtimes \Gamma)^t$ is a finite index inclusion. Sections \ref{sec.IPP} and \ref{sec.bifinite} will imply that we may assume that $\psi(N) \subset N^t$ and that the latter is a finite index inclusion. This allows to prove Theorems \ref{thm.exists} and \ref{thm.main} in Section \ref{sec.proof}. Theorem \ref{thm.exists} follows once we have proven the existence of groups $\Gamma_0,\Gamma_1$ without non-trivial finite-dimensional unitary representations, and actions of these groups on the hyperfinite II$_1$ factor $R$ satisfying all conditions in Theorem \ref{thm.main}. In order to prove this existence, we have to establish in Section~\ref{sec.make-free} the following result: if $\cF_1$ and $\cF_2$ are countable fusion subalgebras of $\FAlg(R)$, where $R$ is the hyperfinite II$_1$ factor, then the set $\al \in \Aut(R)$ such that $\al \cF_1 \al^{-1}$ and $\cF_2$ are free, is a $G_\delta$-dense subset of $\Aut(R)$. This last result generalizes A.3.2 in \cite{IPP}

\section{Results of Ioana, Peterson and Popa and some consequences} \label{sec.IPP}

Throughout this section, we fix infinite groups $\Gamma_0$ and $\Gamma_1$. We set $\Gamma = \Gamma_0 * \Gamma_1$ and take an outer action $\Gamma \actson N$ of $\Gamma$ on the II$_1$ factor $N$. We set $M = N \rtimes \Gamma$, with subalgebras $M_i = N \rtimes \Gamma_i$.

We record from \cite{IPP} the following result. The first statement follows from \cite{IPP}, Theorem 4.3 and the second one from \cite{IPP}, Theorem 1.2.1.

\begin{theorem}[Ioana-Peterson-Popa, \cite{IPP}] \label{thm.IPP}
The following results hold.
\begin{itemize}
\item If $Q \subset M$ is a von Neumann subalgebra with the relative property (T), there exists $i \in \{0,1\}$ such that $Q \embed{M} M_i$.
\item If $t > 0$, $i \in \{0,1\}$ and if $Q \subset M_i^t$ is a von Neumann subalgebra such that $Q \notembed{M_i^t} N^t$, then the quasi-normalizer of $Q$ inside $M^t$ is contained in $M_i^t$.
\end{itemize}
\end{theorem}

\begin{corollary} \label{cor.corIPP}
Suppose that $t > 0$ and that $Q \subset M^t$ is a subfactor with the relative  property (T) whose quasi-normalizer has finite index in $M^t$.
Then, $Q \embed{M^t} N^t$.
\end{corollary}
\begin{proof}
Set $M_i = N \rtimes \Gamma_i$.
Replacing $Q$ by $Q^{1/t}$, we may assume that $t=1$. Suppose that $Q \notembed{M} N$. The first statement in \ref{thm.IPP} yields $i \in \{0,1\}$ such that $Q \embed{M} M_i$. Take a projection $p \in N^n$, a unital $^*$-homomorphism $\psi : Q \recht p M_i^n p$ and a non-zero partial isometry $v \in (\M_{1,n}(\C) \ot M)p$ satisfying $x v = v \psi(x)$ for all $x \in Q$. By construction, the bimodule
$$_{\psi(Q)} \bigl( p (L^2(M_i)^{\oplus n}) \bigr)_{M_i} $$
is isomorphic with a sub-bimodule of $_Q L^2(M)_{M_i}$. Since we are supposing that $Q \notembed{M} N$, we get that $\psi(Q) \notembed{p M_i^n p} p N^n p$. Denote by $Q_1$ the quasi-normalizer of $\psi(Q)$ inside $p M^n p$.
The second statement of Theorem \ref{thm.IPP} implies that $Q_1 \subset p M_i^n p$. But, if $Q_0$ denotes the quasi-normalizer of $Q$ inside $M$, it is clear that $v^* Q_0 v \subset Q_1$. Since we assume that $Q_0$ has finite index in $M$, we arrive at a contradiction.
\end{proof}

The following result is a first step towards the main theorem of the paper.

\begin{theorem} \label{thm.crucialstep}
Let $\Gamma_0$ and $\Gamma_1$ be infinite groups, $\Gamma=\Gamma_0 * \Gamma_1$ their free product and $\Gamma \actson N$ an outer action on the II$_1$ factor $N$. Set $M = N \rtimes \Gamma$ and suppose that $N \subset M$ has the relative property (T).

If $t > 0$ and $\pi : M \recht M^t$ is a finite index, irreducible inclusion, then
$$\pi(N) \embed{M^t} N^t \quad\text{and}\quad N^t \embed{M^t} \pi(N) \; .$$
\end{theorem}
\begin{proof}
By Corollary \ref{cor.corIPP}, we get that $\pi(N) \embed{M^t} N^t$.

Realize $M^t = p M^n p$. Since $\pi(M) \subset M^t$ has finite index, we can take a projection $p_1 \in \pi(M)^m$, a finite index inclusion $\psi : M^t \recht p_1 \pi(M)^m p_1$ and a non-zero partial isometry $v \in p(\M_{n,m}(\C) \ot M)p_1$ satisfying $x v = v \psi(x)$ for all $x \in M^t$. Write $\pi(M)^s :=p_1 \pi(M)^m p_1$. Cutting down if necessary, we may assume that $E_{\pi(M)^s}(v^* v)$ has support $p_1$.

Then, $\psi(N^t) \subset \pi(M)^s$ has the relative property (T). The quasi-normalizer of $\psi(N^t)$ inside $\pi(M)^s$ contains $\psi(M^t)$ and hence, is of finite index. By Corollary \ref{cor.corIPP}, we get that $\psi(N^t) \embed{\pi(M)^s} \pi(N)^s$. So, we find a projection $p_2 \in \pi(N)^k$, a unital $^*$-homomorphism $\theta : \psi(N^t) \recht p_2 \pi(N)^k p_2$ and a non-zero partial isometry $w \in p_1 (\M_{m,k}(\C) \ot \pi(M))p_2$ satisfying $x w = w \theta(x)$ for all $x \in \psi(N^t)$.

Since $E_{\pi(M)^s}(v^*v)$ has support $p_1$ and since $w$ has coefficients in $\pi(M)$, it follows that $vw \neq 0$. Moreover, $N^t vw \subset vw \pi(N)^k$. We have shown that $N^t \embed{M^t} \pi(N)$.
\end{proof}

\section{Bifinite bimodules between crossed products \\ and almost normalizing bimodules} \label{sec.bifinite}

The aim of this section is twofold. First of all, Propositions \ref{prop.propone} and \ref{prop.proptwo} describe the structure of irreducible bifinite $(P \rtimes \Lambda)$-$(N \rtimes \Gamma)$-bimodule \emph{containing a bifinite $P$-$N$-subbimodule.}

The condition of containing a bifinite $P$-$N$-subbimodule is of course a very strong one. Typically, an application of the deformation/rigidity techniques explained in Section \ref{sec.IPP}, yields the existence of a $P$-$N$-subbimodule of finite $N$-dimension and the existence of another $P$-$N$-subbimodule of finite $P$-dimension. In Proposition \ref{prop.assembly}, we show that in good cases this suffices to get the existence of a bifinite $P$-$N$-subbimodule.

Note that Proposition \ref{prop.proptwo} is a generalization of Lemma 8.4 in \cite{IPP}, but we avoid the use of Connes' result about vanishing of $1$-cocycles for finite group actions.

\begin{proposition} \label{prop.propone}
Let $M_0$ be a II$_1$ factor with regular subfactor $N_0$. Suppose that $\Gamma \actson N$ is an outer action of the \emph{ICC group} $\Gamma$ on the II$_1$ factor $N$. Let $H$ be an irreducible bifinite $M_0$-$(N \rtimes \Gamma)$-bimodule containing a bifinite $N_0$-$N$ subbimodule.

Then, there exists a projection $p \in N^n$ and an irreducible finite index inclusion $\psi : M_0 \recht p(N \rtimes \Gamma)^n p$ satisfying
\begin{itemize}
\item $H \cong H(\psi)$ as $M_0$-$(N \rtimes \Gamma)$-bimodules;
\item $\psi(N_0) \subset pN^n p$ and this inclusion has finite index;
\item the relative commutant $p (N \rtimes \Gamma)^n p \cap \psi(N_0)'$ equals $p N^n p \cap \psi(N_0)'$.
\end{itemize}
\end{proposition}

\begin{remark} \label{rem.almostnorm}
The method of the proof below also yields the following result, clarifying the notion of a bifinite bimodule almost normalizing $\Gamma \actson N$. It follows that given such an almost normalizing bifinite $N$-$N$-bimodule $K$, there exists a finite index subgroup $\Gamma_0 < \Gamma$ such that for every $g \in \Gamma_0$, there exist $h,k \in \Gamma$ satisfying
$$H(\si_g) \ot_N K \cong K \ot_N H(\si_h) \quad\text{and}\quad K \ot_N H(\si_g) \cong H(\si_k) \ot_N K \; .$$
See \ref{not.Hpsi} for the notation $H(\si_g)$.
\end{remark}

\begin{proof}[Proof of Proposition \ref{prop.propone}]
Let $H$ be an irreducible bifinite $M_0$-$(N \rtimes \Gamma)$-bimodule containing a bifinite $N_0$-$N$ subbimodule. Since $N \subset N \rtimes \Gamma$ is irreducible, the von Neumann algebra $A$ consisting of $M_0$-$N$-bimodular operators on $H$ is finite-dimensional. Since the elements of $A$ are $M_0$-modular, we write $A$ as acting on the right on $H$.

Take an irreducible bifinite $N_0$-$N$-subbimodule $K \subset H$. Define $\cH$ as the closed linear span of $M_0 K A$. We denote by $z$ the orthogonal projection onto $\cH$ and observe that $z \in \cZ(A)$. Whenever $v \in \cU(A)$, $Kv \cong K$ as $N_0$-$N$-bimodules. So, the regularity of $N_0 \subset M_0$ ensures that $\cH$ is a direct sum of $N_0$-$N$-bimodules isomorphic with one of the $u K$ for $u \in \cN_{M_0}(N_0)$.

Since $\cZ(A)$ is a finite-dimensional abelian algebra normalized by the unitaries $u_g, g \in \Gamma$, we can define the finite index subgroup $\Gamma_0 < \Gamma$ consisting of $g \in \Gamma$ such that $z$ and $u_g$ commute. Hence, for $g \in \Gamma_0$, we have $K u_g \subset \cH$, implying that there exists $u \in \cN_{M_0}(N_0)$ satisfying $Ku_g \cong u K$ as $N_0$-$N$-bimodules. Next define the subset $I \subset \Gamma$ as
$$I:=\{ g \in \Gamma \mid K u_g \cong K \;\;\text{as $N_0$-$N$-bimodules}\; \} \; .$$
It is easily checked that $I$ is globally normalized by the elements of $\Gamma_0$. Moreover, if $g \in I$, we have that $H(\si_g)$ is contained in $\overline{K} \ot_{N_0} K$, implying that $I$ is finite. The ICC property of $\Gamma$ yields that $I = \{e\}$.

Set $M = N \rtimes \Gamma$.
Take an irreducible finite index inclusion $\theta : M_0 \recht qM^m q$ such that $H \cong H(\theta)$ as $M_0$-$M$-bimodules. The presence of $K \subset H$ is then translated to the existence of a non-zero partial isometry $v \in q(\M_{m,n}(\C) \ot M)p_1$ and an irreducible finite index inclusion $\psi_1 : N_0 \recht p_1 N^n p_1$ such that
\begin{itemize}
\item $\theta(x) v = v \psi_1(x)$ for all $x \in N_0$,
\item $K \cong H(\psi_1)$ as $N_0$-$N$-bimodules.
\end{itemize}
We claim that $p_1 M^n p_1 \cap \psi_1(N_0)' = \C p_1$. Indeed, if $\sum_{g \in \Gamma} x_g u_g$ with $x_g \in p_1 N^n \si_g(p_1)$ commutes with $\psi_1(N_0)$, it follows that
$$x_g \si_g(\psi_1(y)) = \psi_1(y) x_g \quad\text{for all}\;\; g \in \Gamma, y \in N_0 \; .$$
So, whenever $x_g \neq 0$, $Ku_g \cong K \ot_N H(\si_g) \cong K$ and hence $g=e$. It follows that our relative commutant lives inside $p_1N^n p_1$ and so, is trivial by the irreducibility of $\psi_1(N_0) \subset p_1 N^n p_1$. The claim is proven.

In particular, we conclude that $v^*v = p_1$ and that $vv^*$ is a minimal projection in $q M^m q \cap \theta(N_0)'$. Also, $v^* \theta(N_0)v \subset p_1 N^n p_1$ and this is a finite index inclusion.

Set $B = qM^m q \cap \theta(N_0)'$. By irreducibility of $\theta(M_0) \subset q M^m q$, we know that $\Ad \theta(\cN_{M_0}(N_0))$ yields an ergodic action on $B$. Since $B$ admits the minimal projection $vv^*$, $B$ is finite-dimensional. Denote by $z$ the central support of $vv^*$ in $B$. Let $(f_{ij})$ be matrix units for $z B$ with $f_{00} = vv^*$. Take a finite set of $u_k \in \cN_{M_0}(N_0)$ such that $\sum_k u_k z u_k^* = q$. Finally, take partial isometries $v_{ki}$ in $N^n$ (enlarging $n$ if necessary) satisfying $v_{ki} v_{ki}^* = p_1$ for all $k,i$ and $p=\sum_{k,i} v_{ki}^* v_{ki}$ a projection in $N^n$. Defining
$$w := \sum_{ki} u_k f_{i0} v v_{ki} \quad\text{and}\quad \psi : M_0 \recht p (N \rtimes \Gamma)^n p : \psi(y) = w^* \theta(y) w \; ,$$
we are done.
\end{proof}

\begin{proposition} \label{prop.proptwo}
Let $\Lambda \overset{\rho}{\actson} P$ and $\Gamma \overset{\sigma}{\actson} N$ be outer actions of the ICC groups $\Lambda,\Gamma$ on the II$_1$ factors $P,N$. Suppose that $H$ is a bifinite $(P \rtimes \Lambda)$-$(N \rtimes \Gamma)$-bimodule containing a bifinite $P$-$N$-subbimodule.

Then there exists an irreducible finite index inclusion $\psi : P \rtimes \Lambda \recht p(N \rtimes \Gamma)^n p$ with $p \in N^n$ and an isomorphism $\delta : \Lambda_0 \recht \Gamma_0$ between finite index subgroups of $\Lambda,\Gamma$, satisfying
\begin{itemize}
\item $H \cong H(\psi)$;
\item $\psi(P) \subset pN^n p$ and this is a finite index inclusion satisfying $p (N \rtimes \Gamma)^n p \cap \psi(P)' = p N^n p \cap \psi(P)'$;
\item for some non-zero projection $z \in \cZ(pN^n p \cap \psi(P)')$, commuting with $\psi(P \rtimes \Lambda_0)$, we have
$$z \psi(u_g) = x_{\delta(g)} u_{\delta(g)} \quad\text{for unitaries}\quad x_s \in z N^n \si_s(z) \;\;\text{when}\;\; s \in \Gamma_0 \; .$$
\end{itemize}
\end{proposition}

\begin{remark}
A close inspection of the proof below, implies that $z$ can be chosen in such a way that the following holds. Define $\psi_0 : P \rtimes \Lambda_0 \recht z (N \rtimes \Gamma_0)^n z : \psi_0(y) = \psi(y)z$ and consider the obvious inclusion bimodules
$$\Hincl(\Lambda_0,\Lambda) = \, _{P \rtimes \Lambda_0} L^2(P \rtimes \Lambda)_{P \rtimes \Lambda} \quad\text{and}\quad
\Hincl(\Gamma_0,\Gamma) = \, _{N \rtimes \Gamma_0} L^2(N \rtimes \Gamma)_{N \rtimes \Gamma} \; .$$
Then, $$H \cong \overline{\Hincl(\Lambda_0,\Lambda)} \underset{P \rtimes \Lambda_0}{\ot} H(\psi_0) \underset{N \rtimes \Gamma_0}{\ot} \Hincl(\Gamma_0,\Gamma) \; .$$
\end{remark}

\begin{proof}[Proof of Proposition \ref{prop.proptwo}]
By Proposition \ref{prop.propone}, we get $H \cong H(\psi)$ where $\psi : P \rtimes \Lambda \recht p (N \rtimes \Gamma)^n p$ is a finite index inclusion satisfying $p \in N^n$, $\psi(P) \subset pN^n p$ a finite index inclusion and $p(N \rtimes \Gamma)^n p \cap \psi(P)' = pN^n p \cap \psi(P)'$.

Let $p_0$ be a minimal projection in the finite dimensional algebra $pN^n p \cap \psi(P)'$ and set $\psi_0(x) = \psi(x)p$ for $x \in P$. Define $K = H(\psi_0)$ as a bifinite $P$-$N$-bimodule. As in the beginning of the proof of Proposition \ref{prop.propone}, we get finite index subgroups $\Lambda_0 < \Lambda$ and $\Gamma_0 < \Gamma$ defined by
\begin{align*}
\Lambda_0 & := \{ g \in \Lambda \mid \; \exists h \in \Gamma \; , \; H(\rho_g) K \cong K H(\si_h) \} \; , \\
\Gamma_0 & := \{ h \in \Gamma \mid \; \exists g \in \Lambda \; , \; K H(\si_h) \cong H(\rho_g) K  \} \; ,
\end{align*}
and an isomorphism $\delta : \Lambda_0 \recht \Gamma_0$ such that $H(\rho_g) K \cong K H(\si_{\delta(g)})$ for all $g \in \Lambda_0$.

Let $z_0 \in \cZ(\psi(P)' \cap pN^n p)$ be the central support of $p_0$. Take $g \in \Lambda_0$. It follows that $\psi(\rho_g(\cdot))z_0$ and $\si_{\delta(g)}(\psi(\cdot) z_0)$ define isomorphic $P$-$N$-bimodules. So, there exists a unitary $v \in \si_{\delta(g)}(z_0)N^n z_0$ such that $v \psi(\rho_g(x)) = \si_{\delta(g)}(\psi(x)) v$ for all $x \in P$. It follows that $u_{\delta(g)}^* v \psi(u_g)$ commutes with $\psi(P)$ and hence, belongs to $pN^n p$. It follows that $z_0 \psi(u_g) \in u_{\delta(g)} N^n$ for all $g \in \Lambda_0$. But then,
$$(\psi(u_h)^* z_0 \psi(u_h)) \; \psi(u_g) = (z_0 \psi(u_h))^* (z_0 \psi(u_{hg}))$$
belongs to $u_{\delta(g)}N^n$ as well, for all $h,g \in \Lambda_0$. Setting $z = \bigvee_{h \in \Lambda_0} \psi(u_h)^* z_0 \psi(u_h)$, we are done.
\end{proof}

The second condition in the next proposition is quite artificial. In the application in this paper, one might as well suppose that $A \subset M$ is a quasi-regular inclusion, i.e.\ $M=\QN_M(A)\dpr$. Elsewhere, we plan another application of the proposition: there it is known that whenever $H \subset L^2(M,\tau)$ is an $A$-$A$-subbimodule with $\dim(H_A) < \infty$, then actually $H \subset L^2(A)$.

\begin{proposition} \label{prop.assembly}
Let $(M,\tau)$ be a von Neumann algebra with faithful normal tracial state $\tau$. Suppose that $A,B \subset M$ are von Neumann subalgebras that satisfy the following conditions.
\begin{itemize}
\item $\dis A \embed{M} B$ and $B \fembed{M} A$.
\item If $H \subset L^2(M,\tau)$ is an $A$-$A$-subbimodule with $\dim(H_A) < \infty$, then $H \subset L^2(\QN_M(A)\dpr)$.
\end{itemize}
Then, there exists a $B$-$A$-subbimodule $K \subset L^2(M,\tau)$ satisfying
$$\dim( _B K) < \infty \quad\text{and}\quad \dim(K_A) < \infty \; .$$
So, there exists a projection $p \in \M_n(\C) \ot A$, a non-zero partial isometry $v \in (\M_{1,n}(\C) \ot M)p$ and a unital $^*$-homomorphism $\theta : B \recht p A^n p$ satisfying
$$\theta(B) \subset pA^n p \;\;\text{has finite index, and}\quad b v = v \theta(b) \quad\text{for all}\quad b \in B \; .$$
\end{proposition}

In the above statement, all dimensions are with respect to the restriction of $\tau$ to $A$ and $B$. In particular, the index of $\theta(B) \subset p A^n p$, is defined as $\dim( L^2(pA^n p)_B)$, where the right $B$-module action is through $\theta$.

\begin{proof}
Denote by $J$ the anti-unitary operator on $L^2(M,\tau)$ given by $J x = x^*$. Then, $J (\la M,e_A \ra \cap B') J = \la M,e_B \ra \cap A'$. So, we get two normal faithful traces on $\la M,e_A \ra \cap B'$: one denoted by $\Tr_A$ and defined by restricting the trace on $\la M,e_A \ra$ and the other denoted by $\Tr_B$ and obtained by applying the previous formula and restricting the trace on $\la M,e_B \ra$. Define
\begin{align}
p &= \vee \{ p_0 \mid p_0 \;\;\text{projection in}\;\; \la M,e_A \ra \cap B' \;\;\text{with}\;\; \Tr_A(p_0) < \infty \} \; , \label{eq.p}\\
q &= \vee \{ q_0 \mid q_0 \;\;\text{projection in}\;\; \la M,e_A \ra \cap B' \;\;\text{with}\;\; \Tr_B(q_0) < \infty \} \; . \notag
\end{align}
It suffices to prove that $pq \neq 0$. Indeed, approximating $p$ and $q$, we get $p_0$ with $\Tr_A(p_0) < \infty$ and $q_0$ with $\Tr_B(q_0) < \infty$, satisfying $p_0 q_0 \neq 0$. Taking a spectral projection of the positive operator $q_0 p_0 q_0$, we arrive at an orthogonal projection $r \in \la M,e_A \ra \cap B'$ satisfying $\Tr_A(r),\Tr_B(r) < \infty$. Taking $K= r L^2(M,\tau)$, the lemma is proved.

Take non-zero partial isometries $v,w \in \M_{1,n}(\C) \ot M$ and, possibly non-unital, $^*$-homomorphisms $\rho : A \recht B^n$, $\theta : B \recht A^n$ such that
$$a v = v \rho(a) \; , \quad b w = w \theta(b) \quad\text{for all}\quad a \in A, b \in B \; .$$
Since $B \fembed{M} A$, we may assume that $v (1 \ot w) \neq 0$. Note that $w w^* \in M \cap B'$, so that we may assume that $v = v(1 \ot ww^*)$.
By construction, the right $A$-module generated by the (finitely many) coefficients of $v(1 \ot w)$, is also a left $A$-module. Our assumptions imply that the coefficients of $v(1 \ot w)$ belong to $\QN_M(A)\dpr$. With $p$ defined by \eqref{eq.p}, it is easily checked that $H_0:=p L^2(M,\tau)$ is a right $\QN_M(A)\dpr$-module. By construction, the coefficients of $w$ belong to $H_0$ and hence, the coefficients of $v^* = w (v (1 \ot w))^*$ belong to $H_0$ as well. By construction, the coefficients of $v^*$ belong to $q L^2(M,\tau)$. So, we have shown that $pq \neq 0$.
\end{proof}

\section{Proof of the main theorems} \label{sec.proof}

\begin{proof}[Proof of Theorem \ref{thm.main}]
Write $\Gamma=\Gamma_0 * \Gamma_1$ and $M = N \rtimes \Gamma$. Let $H$ be a bifinite $M$-$M$-bimodule.
Combining Theorem \ref{thm.crucialstep}, Proposition \ref{prop.assembly} and Proposition \ref{prop.proptwo}, we get $H \cong H(\psi)$ where $\psi : M \recht p M^n p$ is an irreducible finite index inclusion satisfying
\begin{itemize}
\item $p \in N^n$ and $\psi(N) \subset pN^n p$ a finite index inclusion,
\item $p M^n p \cap \psi(N)' = pN^n p \cap \psi(N)'$,
\item $\psi(u_g) z = x_{\delta(g)} u_{\delta(g)}$ for all $g \in \Lambda$, where $\Lambda < \Gamma$ is a finite index subgroup, $\delta : \Lambda \recht \Gamma$ an injective homomorphism with finite index image, $x_h$ a unitary in $z N^n \si_h(z)$ for all $h \in \delta(\Lambda)$ and $z$ a central projection in $pN^n p \cap \psi(N)'$ commuting with $\psi(N \rtimes \Lambda)$.
\end{itemize}
Denote by $K$ the bifinite $N$-$N$-bimodule defined by the inclusion $N \recht zN^n z : x \mapsto \psi(x) z$. We prove that $K$ is a multiple the trivial $N$-$N$-bimodule, which will almost end the proof of the theorem.

Set $\Lambda_i := \Gamma_i \cap \Lambda$ and note that $\Lambda_i$ is a finite index subgroup of $\Gamma_i$. We assumed that $\Gamma_0,\Gamma_1,\Z$ have no isomorphic finite index subgroups and that the finite index subgroups of $\Gamma_0,\Gamma_1$ are freely indecomposable. Hence, the Kurosh theorem implies that $\delta(\Lambda_i)$ is a finite index subgroup of $s_i \Gamma_i s_i^{-1}$, for some $s_0,s_1 \in \Gamma$. Unitary conjugating with $u_{s_0}$ from the beginning, we may assume that $\delta(\Lambda_0)$ is a finite index subgroup of $\Gamma_0$ and that $\delta(\Lambda_1)$ is a finite index subgroup of $s \Gamma_1 s^{-1}$. Again unitary conjugating, we may assume that either $s=e$ or $s \in (\Gamma_1 - \{e\}) \cdots (\Gamma_0 - \{e\})$.

So, the map $N \rtimes \Lambda_0 \recht z(N \rtimes \Lambda_0)^n z : y \mapsto \psi(y)z$ defines a bifinite $(N \rtimes \Lambda_0)$-$(N \rtimes \Lambda_0)$-bimodule that contains the $N$-$N$-bimodule $K$. By Lemma \ref{lem.almostnorm}, $K$ is almost normalizing $\Gamma_0 \actson N$. By our assumptions $K \cup \Gamma_0$ and $\Gamma_1$ are free inside $\FAlg(N)$. Writing for all $g \in \Lambda_1$, $\delta(g) = s \eta(g) s^{-1}$ for $\eta(g) \in \Gamma_1$ and $s$ as above, the formula $\psi(u_g) z = x_{\delta(g)} u_{\delta(g)}$ implies that $H(\si_g) K \cong K H(\si_{s\eta(g)s^{-1}})$ for every $g \in \Lambda_1$. Given the form of $s$, this is a contradiction with the freeness of $K \cup \Gamma_0$ and $\Gamma_1$, unless $K$ is a multiple of the trivial $N$-$N$-bimodule.

Our claim is proven and we find a non-zero partial isometry $v \in p(\M_{n,1}(\C) \ot N)$ satisfying
\begin{equation}\label{eq.formuleke}
\psi(x) v = vx \;\;\text{for all}\;\; x \in N \; .
\end{equation}
Then, $v^*v = 1$ and \eqref{eq.formuleke} remains true replacing $v$ by $q \psi(u_g) v u_g^*$ whenever $g \in \Gamma$ and $q \in pN^n p \cap \psi(N)'$. It follows that we can find a unitary $w \in p(\M_{n,k}(\C) \ot N)$ satisfying $w^* \psi(x) w = 1 \ot x$ for all $x \in N$. It is now an exercise to check that $w^* \psi(u_g) w = \theta(g) \ot u_g$ for some representation $\theta : \Gamma \recht \cU(k)$.
\end{proof}

Finally, we prove the existence of groups and actions satisfying all the requirements in Theorem \ref{thm.main} and moreover such that the groups do not admit finite-dimensional unitary representations.

\begin{proof}[Proof of Theorem \ref{thm.exists}]
We have to prove that there exist infinite groups $\Gamma_0,\Gamma_1$ together with outer actions on the hyperfinite II$_1$ factor $R$ such that all conditions of Theorem \ref{thm.main} are satisfied and such that all finite dimensional unitary representations of $\Gamma_0$ and $\Gamma_1$ are trivial.

Consider the group $A_\infty$ of finite even permutations of $\N$. It is well known that every finite dimensional unitary representation of $A_\infty$ is trivial. Consider $\Z/3\Z \subset A_\infty$, identifying $1$ and the cyclic permutation of $\{0,1,2\}$. Finally, consider $\Z/3\Z \subset \SL(3,\Z)$ identifying $1$ and the matrix $\Bigl(\begin{smallmatrix} 0 & 0 & 1 \\ 1 & 0 & 0 \\ 0 & 1 & 0 \end{smallmatrix}\Bigr)$. We then define
$$\Gamma_0 = \SL(3,\Z) \underset{\Z/3\Z}{*} A_\infty \quad\text{and}\quad \Gamma_1 = A_\infty \; .$$
As stated above, $\Gamma_1$ does not have non-trivial finite dimensional unitary representations. If $\pi : \Gamma_0 \recht \cU(n)$ is a finite dimensional unitary representation, $A_\infty \subset \Ker \pi$. In particular, $\Z/3\Z \subset \Ker \pi$. Since the smallest normal subgroup of $\SL(3,\Z)$ containing $\Z/3\Z$, is the whole of $\SL(3,\Z)$, it follows that $\Ker \pi = \Gamma_0$.

In particular, $\Gamma_0$ and $\Gamma_1$ do not have non-trivial finite index subgroups. Both $\SL(3,\Z)$ and $A_\infty$ are freely indecomposable. Then, the Kurosh theorem implies that $\Gamma_0$ is freely indecomposable as well.

We next claim that there exists an outer action of $\Gamma_0$ on the hyperfinite II$_1$ factor $R$ such that $R \subset R \rtimes \Gamma_0$ has the relative property (T). First take an outer action of $\SL(3,\Z)$ on $R$ such that $R \subset R \rtimes \SL(3,\Z)$ has the relative property (T). A way of doing so, goes as follows. Consider the semi-direct product $\SL(3,\Z) \ltimes (\Z^3 \times \Z^3)$, where $A \cdot (x,y) = (Ax,(A^{-1})^t y)$ for all $A \in \SL(3,\Z)$ and $x,y \in \Z^3$. It is clear that $\Z^3 \times \Z^3$ is a subgroup with the relative property (T). Take an $\SL(3,\Z)$-invariant non-degenerate $2$-cocycle $\om$ on $\Z^3 \times \Z^3$. We then get the required action of $\SL(3,\Z)$ on $R = \cL_\om(\Z^3 \times \Z^3)$. Next, take any outer action of $A_\infty$ on $R$. By Connes' uniqueness theorem for outer actions of finite cyclic groups on $R$ (see \cite{C2}), we may assume that the actions of $\Z/3 \Z \subset A_\infty$ and $\Z / 3\Z \subset \SL(3,\Z)$ coincide. Hence, we get
an action of $\Gamma_0$ on $R$. Further modifying the action of $A_\infty$ by applying Proposition \ref{prop.free-amalgamation}, we have shown that there exists an outer action of $\Gamma_0$ on $R$ that extends the $\SL(3,\Z)$ action. Then, $R \subset R \rtimes \Gamma_0$ still has the relative property (T).

Finally, take any outer action of $\Gamma_1$ on the hyperfinite II$_1$ factor $R$. Denote by $\cF$ the fusion subalgebra of $\FAlg(R)$ generated by the bifinite $R$-$R$-bimodules almost normalizing $\Gamma_0 \actson R$. By Lemma \ref{lem.countability} below, $\cF$ is countable. It follows from Theorem \ref{thm.make-free} below that there exists an automorphism $\al \in \Aut(R)$ such that $\cF$ and $\al \Gamma_1 \al^{-1}$ are free in the sense of Definition \ref{def.free}. Replacing $\Gamma_1$ by $\al \Gamma_1 \al^{-1}$, all conditions of Theorem \ref{thm.main} are fulfilled and moreover, $\Gamma$ only has trivial finite dimensional unitary representations. So, we are done.
\end{proof}

\begin{lemma}\label{lem.countability}
Let $N$ be a II$_1$ factor and $\Gamma \actson N$ an outer action such that
$N \subset N \rtimes \Gamma$ has the relative property (T). Then, the almost normalizer of $\Gamma \actson N$ in $\FAlg(N)$ (in the sense of Definition \ref{def.almostnorm}) is countable.
\end{lemma}
\begin{proof}
Set $M = N \rtimes \Gamma$.
By contradiction and countability of $\Gamma$ and $\N$, it is sufficient to prove the following:
if $n \in \N_0$ and if $\psi_i : M \recht p_i M^n p_i$ defines an uncountable family of bifinite $M$-$M$-bimodules $H_i$ containing non-zero irreducible bifinite $N$-$N$-bimodule $K_i \subset H_i$, there exist $i \neq j$ and $g,h \in \Gamma$ such that $K_i \cong H(\si_g) K_j H(\si_h)$ as $N$-$N$-bimodules.

Take $\eps > 0$ and $F \subset M$ finite such that every $M$-$M$-bimodule $H$ that admits a vector $\xi \in H$ with the properties $1-\eps \leq \|\xi\| \leq 1$ and $|\langle \xi , a \xi b\rangle - \tau(ab)| < \eps$ for all $a,b \in F$, actually admits a non-zero $N$-central vector.

Assume for convenience that $1 \in F$ and consider the $\psi_i$ as non-unital homomorphisms $M \recht M^n$. By the pigeon hole principle, we can find $i \neq j$ such that $\|\psi_i(x) - \psi_j(x)\|_2 < \eps \|q_i\|_2$ for all $x \in F$. Consider the $M$-$M$-bimodule $p_i L^2(M^n) p_j$ with left action given by $\psi_i$ and right action by $\psi_j$. The vector $\xi = \|p_i\|_2^{-1} p_i p_j$ satisfies the above conditions and we conclude that $p_i L^2(M^n) p_j$ contains a non-zero $N$-central vector. It follows that there exist irreducible $N$-$N$-subbimodules $\Ktil_i \subset H_i$ and $\Ktil_j \subset H_j$ with $\Ktil_i \cong \Ktil_j$ as $N$-$N$-bimodules. To conclude to proof, it suffices to observe that for every $i$, $H_i$ as an $N$-$N$-bimodule is a direct sum of irreducible $N$-$N$-bimodules isomorphic with $H(\si_g) K_j H(\si_h)$, $g, h \in \Gamma$.
\end{proof}

\section{Realizing freeness between fusion subalgebras of $\FAlg(R)$} \label{sec.make-free}

In this section, we prove the following crucial result: whenever $\cF_1,\cF_2$ are \emph{countable} fusion subalgebras of $\FAlg(R)$, where $R$ denotes the hyperfinite II$_1$ factor, there exists an automorphism $\al \in \Aut(R)$ such that $$\cF_1^\al := H(\al^{-1})\cF_1 H(\al) \quad\text{and}\quad \cF_2$$ are free. In the terminology of \cite{BJ}, this implies that any two hyperfinite finite index subfactors admit a hyperfinite realization of their free composition (see page 94 in \cite{BJ}).

\begin{theorem} \label{thm.make-free}
Let $R$ be the hyperfinite II$_1$ factor.
Let $\cF_1,\cF_2$ be countable fusion algebras of bifinite $R$-$R$-bimodules. Then,
$$\{\al \in \Aut(R) \mid \cF_1^\al \;\;\text{and}\;\; \cF_2 \;\; \text{are free}\}$$
is a $G_\delta$ dense subset of $\Aut(R)$.
\end{theorem}

Recall that if $_M H_M$ is an $M$-$M$-bimodule and $A \subset M$ a von Neumann subalgebra, a vector $\xi \in H$ is said to be \emph{$A$-central} if $a \xi = \xi a$ for all $a \in A$. Note that if $p$ denotes the orthogonal projection onto the subspace of $A$-central vectors, $p \xi$ is precisely the element of minimal norm in the closed convex hull
$$\overline{\operatorname{co}} \{ u \xi u^* \mid u \in \cU(A) \} \; .$$

In what follows, we make use of the following special property for a bifinite bimodules $_R H_R$ over the \emph{hyperfinite} II$_1$ factor $R$. Fix a free ultrafilter $\om$ on $\N$ and consider the ultrapower algebra $R^\om$. We claim that there exists $n \in \N$ and an $R$-$R$-bimodular isometric embedding $V : H \recht L^2(R^\om)^{\oplus n}$ into the $n$-fold direct sum of $_R L^2(R^\om)_R$. Denoting by $\cH$ the W$^*$-bimodule of bounded vectors in $H$, we can take $V \cH \subset \M_{n,1}(\C) \ot R^\om$. To prove the existence of such an embedding, take $\psi : R \recht p R^n p$ such that $H \cong H(\psi)$. We can take a partial isometry $A \in \M_n(\C) \ot R^\om$ satisfying $A^* A = p$ and $(1 \ot x) A = A \psi(x)$ for all $x \in R$. It then suffices to define
$$p(L^2(R)^{\oplus n}) \recht L^2(R^\om)^{\oplus n} : \xi \mapsto A \xi \; .$$
Moreover, $_R H_R$ does not contain the trivial bimodule if and only if $(\id \ot E_{R' \cap R^\om})(V \xi) = 0$ for all $\xi \in \cH$.

We are now ready to prove Theorem \ref{thm.make-free} and the proof will be based on the technical Proposition \ref{prop.unitary} below.

\begin{proof}[Proof of Theorem \ref{thm.make-free}]
Suppose that $H_0,\ldots,H_{2k}$ are irreducible bifinite $R$-$R$-bimodules, with $H_j$ non-trivial if $1 \leq j \leq 2k-1$.
When $\al \in \Aut(R)$ and $H \in \FAlg(R)$, we write $H^\al := H(\al^{-1}) H H(\al)$ and define
$$K(\al):= H_0 \, H_1^\al \, H_2 \, H_3^\al \, \cdots \, H_{2k-1}^\al \, H_{2k} \; .$$
We have to prove that
$$W:=\{\al \in \Aut(R) \mid K(\al) \;\;\text{is disjoint from the trivial bimodule}\}$$
is a $G_\delta$ dense subset of $\Aut(R)$.

Let $\cH_i \subset H_i$ denote the W$^*$-$M$-$M$-bimodules that sit densely inside $H_i$. Take $n$ sufficiently large and take isometric embeddings
$$V_i : H_i \recht L^2(R^\om)^{\oplus n} \quad\text{with}\quad V_i \cH_i \subset \M_{n,1}(\C) \ot R^\om \; .$$
Denote by $\proj$ the orthogonal projection onto the $R$-central vectors of $_R K(\al)_R$. Whenever $\xi_i \in \cH_i$ and $\eps > 0$, we define
$$W(\xi_0,\ldots,\xi_{2k} \;  ; \; \eps) := \{\al \in \Aut(R) \mid \|\proj(\xi_0 \ot \cdots \ot \xi_{2k})\| < \eps \} \; .$$
We prove three statements.
\begin{enumerate}
\item Every $W(\xi_0,\ldots,\xi_{2k} \; ; \; \eps)$ is open in $\Aut(R)$.
\item Every $W(\xi_0,\ldots,\xi_{2k} \; ; \; \eps)$ is dense in $\Aut(R)$.
\item Taking the intersection of $W(\xi_0,\ldots,\xi_{2k} \; ; \; \frac{1}{m})$ where $m$ runs through $\N_0$ and the $\xi_i$ run through a countable $\|\cdot\|_2$-dense subset of $\cH_i$, we precisely obtain $W$.
\end{enumerate}
By the Baire category theorem, these statements together show that $W$ is a $G_\delta$ dense subset of $\Aut(R)$.

To prove the first statement, observe that $W(\xi_0,\ldots,\xi_{2k} \; ; \; \eps)$ is the union of all
$$\Bigl\{\al \in \Aut(R) \; \Big| \; \Bigl\| \sum_{i=1}^n \lambda_i u_i (\xi_0 \ot \cdots \ot \xi_{2k}) u_i^* \Bigr\|_{K(\al)} < \eps \Bigr\} \; ,$$
where $n$ runs through $\N_0$, where $\lambda_1,\ldots,\lambda_n$ runs through all $n$-tuples of positive real numbers with sum $1$ and where $u_1,\ldots,u_n$ runs through all $n$-tuples of unitaries in $R$. All these sets are easily seen to be open.

To prove the second statement, set $V_i \xi_i = y_i = (y_i(1),\ldots,y_i(n))^t \in \M_{n,1}(\C) \ot R^\om$. Then, extending an automorphism of $R$ to an automorphism of $R^\om$ in the canonical way, we have
\begin{equation}\label{eq.joepie}
\|\proj(\xi_1 \ot \cdots \ot \xi_{2k})\|^2 = \sum_{i_0,\ldots,i_{2k}=1}^n \|E_{R' \cap R^\om}(\, y_0(i_0) \; \al(y_1(i_1)) \; y_2(i_2) \; \cdots \; \al(y_{2k-1}(i_{2k-1})) \; y_{2k}\,)\|_2^2 \; .
\end{equation}
Fix $\be \in \Aut(R)$. We show that $\be$ is in the closure of
$W(\xi_0,\ldots,\xi_{2k} \; ; \; \eps)$. Write $R$ as the infinite tensor product of $2$ by $2$ matrices, yielding $R = \M_{2^s}(\C) \ot R_s$. It is sufficient to prove that, for every $s \in \N_0$, there exists a unitary $u \in R_s$ such that $(\Ad u) \be \in W(\xi_0,\ldots,\xi_{2k} \; ; \; \eps)$. The existence of $u$ follows combining \eqref{eq.joepie}, Proposition \ref{prop.unitary} and the following observations.
\begin{itemize}
\item If $H_i$ is disjoint from the trivial bimodule and $\be \in \Aut(R)$ is arbitrary, $H^\be_i$ does not admit non-zero $R$-central vectors either and hence, does not even admit non-zero $R_s$-central vectors. So, $$E_{R_s' \cap R^\om}(\be(y_i(j))) = 0$$ for all $j = 1,\ldots,n$, all $s$ and all $\be \in \Aut(R)$.
\item By construction, the elements $\be(y_i(j)) \in R^\om$ quasi-normalize $R$ for all $\be \in \Aut(R)$. Hence, they quasi-normalize $R_s$ for all $s$.
\item We apply Proposition \ref{prop.unitary} to the subfactor $R_s$ of the von Neumann algebra generated by $R$, the $y_{2i}(j)$ and $\be(y_{2i+1}(j))$.
\end{itemize}

It remains to prove the third statement. Of course, if $\al \in W$, then $\al \in W(\xi_0,\ldots,\xi_{2k} \; ; \; \eps)$ for all $\xi_i$ and $\eps > 0$.  Conversely, if $\al$ belongs to the intersection stated above, we have
$$\proj(\xi_0 \ot \cdots \ot \xi_{2k}) = 0$$
for dense families of $\xi_i \in H_i$. But this implies that $\proj=0$ and so $\al \in W$.
\end{proof}

We have the following variant of Theorem \ref{thm.make-free}, that we use in the proof of Theorem \ref{thm.exists}.

\begin{proposition} \label{prop.free-amalgamation}
Suppose that the countable groups $\Gamma_0,\Gamma_1$ have a common finite subgroup $K$. Let $\Gamma_0 \underset{K}{*} \Gamma_1$ act on the hyperfinite II$_1$ factor $R$. Suppose that both $\Gamma_0$ and $\Gamma_1$ act outerly. Denote by $\Aut_K(R)$ the automorphisms of $R$ that commute with all the automorphisms in $K$. Then,
$$\{ \al \in \Aut_K(R) \mid \; \text{The subgroups}\;\;\Gamma_0 \;\;\text{and}\;\; \al \Gamma_1 \al^{-1} \;\;\text{of} \;\; \Out(R)\;\; \text{are free with amalgamation over}\;\; K \}$$
is a $G_\delta$ dense subset of $\Aut_K(R)$.
\end{proposition}
\begin{proof}
One can almost entirely copy the proof of Theorem \ref{thm.make-free}, using the following observation. Let $\al \in \Aut(R)$ be such that $\si_k \al$ is outer for every $k \in K$. Denote by $R^K$ the fixed point algebra of $K$. We claim that the $R$-$R$-bimodule $H(\al)$ has no non-zero $R^K$-central vectors. If it would, the irreducibility of $R^K \subset R$ implies that there exists a unitary $v \in R$ such that $v\al(x)v^* = x$ for all $x \in R^K$. By Jones' uniqueness theorem for outer actions of finite groups (see \cite{Jon1}), we may assume that the action of $K$ is dual and conclude that $(\Ad v) \al = \si_k$ for some $k \in K$. This contradicts our assumption and proves that $H(\al)$ has no non-zero $R^K$-central vectors. Writing $R^K$ as an infinite tensor product of $2$ by $2$ matrices, we get $R^K = M_{2^k}(\C) \ot R_k$.
If $A \in R^\om$ is a unitary implementing $\al$, it follows as in the proof of \ref{thm.make-free} that $E_{R_k' \cap R^\om}(A) = 0$. This is again the starting point to apply Proposition \ref{prop.unitary}.
\end{proof}

The following is the crucial result to obtain Theorem \ref{thm.make-free}. Most of the proof is taken almost literally from Lemmas 1.2, 1.3 and 1.4 in \cite{P10}. We repeat the argument for the convenience of the reader, since slight modifications are needed: in \cite{P10}, the relative commutant $N' \cap M$ is assumed to be finite-dimensional, while we assume that $N$ is a factor and the inclusion $N \subset M$ quasi-regular. This forces us to prove the extra lemma \ref{lem.eigen} below.

\begin{proposition} \label{prop.unitary}
Let $(M,\tau)$ be a von Neumann algebra with faithful normal tracial state $\tau$. Let $N \subset M$ be a von Neumann subalgebra. Suppose that $N$ is a factor of type II$_1$ and that $N$ is quasi-regular in $M$. Let $\cV \subset M$ be a finite subset such that $E_{N' \cap M}(A) = 0$ for all $A \in \cV$.

For every $\eps > 0$ and every $K \in \N_0$, there exists a unitary $u \in N$ such that
$$\|E_{N' \cap M}(A_0 u^{k_1} A_1 u^{k_2} A_2 \cdots u^{k_n} A_n) \|_2 < \eps$$
for all $1 \leq n \leq K$, $1 \leq |k_i| \leq K$, $A_0,A_n \in \cV \cup \{1\}$ and $A_1,\ldots,A_{n-1} \in \cV$.
\end{proposition}

Proposition \ref{prop.unitary} is proven below, after the following preliminary result.

\begin{lemma} \label{lem.partial}
Let $(M,\tau)$ be a von Neumann algebra with faithful normal tracial state $\tau$. Let $N \subset M$ be a von Neumann subalgebra. Suppose that $N$ is a factor of type II$_1$ and that $N$ is quasi-regular in $M$. Let $f \in N$ be a non-zero projection and $\cV \subset M$ a finite subset such that $E_{N' \cap M}(fAf) = 0$ for all $A \in \cV$.

For every $\eps > 0$ and every $K \in \N_0$, there exists a partial isometry $v \in fNf$ satisfying $vv^* = v^* v$, $\tau(vv^*) \geq \tau(f)/4$ and
$$\|E_{N' \cap M}(A_0 v^{k_1} A_1 v^{k_2} A_2 \cdots v^{k_n} A_n) \|_2 < \eps$$
for all $1 \leq n \leq K$, $1 \leq |k_i| \leq K$, $A_0,A_n \in \cV \cup \{1\}$ and $A_1,\ldots,A_{n-1} \in \cV$.
\end{lemma}

Here, and in what follows, we use the convention that $v^0 := vv^*$ and $v^{-k} := (v^*)^k$ for $k \in \N_0$, whenever $v$ is a partial isometry satisfying $vv^* = v^* v$.

\begin{proof}
We may assume that $\|A\| \leq 1$ for all $A \in \cV$. Since $\|z\|_2^2 \leq \|z\| \, \|z\|_1$, we prove the following:
for every $\eps > 0$ and every $K \in \N_0$, there exists a partial isometry $v \in fNf$ such that $vv^* = v^* v$, $\tau(vv^*) \geq \tau(f)/4$ and
$$\|E_{N' \cap M}(A_0 v^{k_1} A_1 v^{k_2} A_2 \cdots v^{k_n} A_n) \|_1 \leq \eps$$
for all $1 \leq n \leq K$, $1 \leq |k_i| \leq K$, $A_0,A_n \in \cV \cup \{1\}$ and $A_1,\ldots,A_{n-1} \in \cV$.

Fix $\eps > 0$ and $K \in \N_0$.
Let $\eps_0 > 0$ and define $\eps_{n} = 2^{n+1} \eps_{n-1}$, up to $\eps_K$. Take $\eps_0$ small enough such that $\eps_K < \eps$. Define $I$ as the set of partial isometries $v \in fNf$ satisfying $vv^* = v^*v$ and
$$
\|E_{N' \cap M}(A_0 v^{k_1} A_1 v^{k_2} A_2 \cdots v^{k_n} A_n) \|_1 \leq \eps_n \tau(vv^*)$$
for all $1 \leq n \leq K$, $1 \leq |k_i| \leq K$, $A_1,\ldots,A_{n-1} \in \cV$, $A_0 \in \cV \cup f\cV \cup \{1\}$ and $A_n \in \cV \cup \cV f \cup \{1\}$.

Order $I$ by inclusion of partial isometries. By Zorn's lemma, take a maximal element $v \in I$ and set $p = v v^*$. It might be that $v=0$. If $\tau(p) \geq \tau(f)/4$, we are done. Otherwise $\tau(p)<\tau(f)/4$ and we set $p_1 := f-p$. Note that $\tau(p)/\tau(p_1) \leq 1/3$. Write $M_1:= p_1 M p_1$, with normalized tracial state $\tau_1$ and corresponding norms $\|\cdot\|_{1,M_1}$ and $\|\cdot\|_{2,M_1}$. Applying Theorem A.1.4 in \cite{P9} to the inclusion $p_1 N p_1 \subset p_1 M p_1$, take a non-zero projection $q \in p_1 N p_1$, such that
$$\|q x q - E_{(N' \cap M)p_1}( p_1 x p_1) q \|_{2,M_1} \leq \eps_0 \|q\|_{2,M_1}$$
for all $x = A_1 v^{k_1} \cdots v^{k_{s-1}} A_s$ and all $1 \leq s \leq K$, $1 \leq |k_i| \leq K$ and $A_1,\ldots,A_s \in \cV$. We shall prove that a unitary $w \in qNq$ can be chosen in such a way that $v+w \in I$. This then contradicts the maximality of $v$.

Let $x = A_1 v^{k_1} \cdots v^{k_{s-1}} A_s$ with $1 \leq s \leq K$, $1 \leq |k_i| \leq K$ and $A_1,\ldots,A_s \in \cV$. Observe that
$$\|q x q - E_{(N' \cap M)p_1}( p_1 x p_1 ) q \|_{1,M_1} \leq \|q x q - E_{(N' \cap M)p_1}( p_1 x p_1 ) q \|_{2,M_1} \|q\|_{2,M_1} \leq \eps_0 \tau_1(q) \; .$$
One checks that $\|E_{(N' \cap M)p_1}( p_1 x p_1 ) q \|_{1,M_1} = \|E_{N' \cap M}(x p_1)\|_1 \; \tau_1(q)/\tau(p_1)$. On the other hand,
\begin{align*}
\|E_{N' \cap M}(x p_1)\|_1 & \leq \|E_{N' \cap M}(xf)\|_1 + \|E_{N' \cap M}(x p)\|_1 = \|E_{N' \cap M}(xf)\|_1 + \|E_{N' \cap M}(v x v^*)\|_1
\\ & \leq (\eps_{s-1} + \eps_{s+1})\tau(p) \; .
\end{align*}
It follows that $\|E_{(N' \cap M)p_1}( p_1 x p_1 ) q \|_{1,M_1} \leq \tau_1(q) \; (\eps_{s-1} + \eps_{s+1})/3$.
Altogether, we conclude that
\begin{equation}\label{eq.finalest}
\|qxq \|_1 \leq \eps_{s+1} \tau(q)/2 \; .
\end{equation}

By Lemma \ref{lem.eigen} below, take a unitary $w \in qNq$ such that
$$
\|E_{N' \cap M}(A_0 v^{k_1} \cdots A_{j-1} w^{k_j} A_j \cdots v^{k_n} A_n)\|_1 \leq \frac{\eps_n \tau(q)}{4n}$$
for all $1 \leq n \leq K$, $1 \leq j \leq n$, $1 \leq |k_i| \leq K$, $A_1,\ldots,A_{n-1} \in \cV$, $A_0 \in \cV \cup f\cV \cup \{1\}$ and $A_n \in \cV \cup \cV f \cup \{1\}$.

{\bf Claim:} the partial isometry $v+w$ belongs to $I$, contradicting the maximality of $v$. To prove the claim, take $1 \leq n \leq K$, $1 \leq |k_i| \leq K$, $A_1,\ldots,A_{n-1} \in \cV$, $A_0 \in \cV \cup f\cV \cup \{1\}$ and $A_n \in \cV \cup \cV f \cup \{1\}$. We develop the sums in the expression
\begin{equation}\label{eq.tussen}
\begin{split}
E_{N' \cap M}(A_0 (v^{k_1}+w^{k_1}) A_1 (v^{k_2}+w^{k_2}) A_2 \cdots (v^{k_n}+w^{k_n}) A_n) \; .
\end{split}
\end{equation}
\begin{itemize}
\item There is one term with only $v$'s appearing. Its $\|\cdot\|_1$-norm is bounded by $\eps_n \tau(p)$, because $v \in I$.
\item There are $n$ terms with $w$ appearing at one place. Each term has its $\|\cdot\|_1$-norm bounded by $\frac{\eps_n \tau(q)}{4n}$. Altogether, their $\|\cdot\|_1$-norm is bounded by $\eps_n \tau(q)/4$.
\item There is $1$ term with $w$ appearing in position $1$ and position $n$ and with $v$'s in the other positions. This term contains the subexpression
$$q A_1 v^{k_2} \cdots v^{k_{n-1}} A_{n-1} q \; .$$
Because of \eqref{eq.finalest}, the $\|\cdot\|_1$-norm of this term is bounded by $\eps_n \tau(q)/2$.
\item There are less than $2^n$ terms where $w$ appears on at least two positions that are not exactly the positions $1,n$. In every such term, we have the subexpression
$$q A_i v^{k_{i+1}} \cdots v^{k_j} A_j q$$
with $1 \leq i \leq j \leq n-1$ and $0 \leq j-i \leq n-3$. By \eqref{eq.finalest}, the $\|\cdot\|_1$-norm of this subexpression is bounded by $\eps_{n-1} \tau(q)/2$. It follows that the sum of all the terms of this type has $\|\cdot\|_1$-norm bounded by $2^{n-1} \eps_{n-1} \tau(q) \leq \eps_n \tau(q)/4$.
\end{itemize}
It follows that the $\|\cdot\|_1$-norm of the expression in \eqref{eq.tussen} is bounded by $\eps_n(\tau(p) + \tau(q)) = \eps_n \tau(p+q)$, proving that $v+w \in I$.
\end{proof}

\begin{lemma}\label{lem.eigen}
Let $(M,\tau)$ be a von Neumann algebra with faithful normal tracial state. Let $N \subset M$ be a von Neumann subalgebra. Suppose that $N$ is a factor of type II$_1$ and that $N$ is quasi-regular in $M$. If $w_n$ is a bounded sequence in $N$ that converges weakly to $0$, then
$$\|E_{N' \cap M}(a w_n b) \|_2 \recht 0$$
for all $a,b \in M$.
\end{lemma}

\begin{proof}
{\bf Step 1.} Let $a \in M$ with $\|a\|\leq 1$. The sequence $\|E_{N' \cap M}(a w_n)\|_2$ converges to $0$, whenever $w_n$ is a bounded sequence in $N$ that converges weakly to zero. Indeed, writing $E_{N' \cap M} = E_{N' \cap M} \circ E_{N \vee (N' \cap M)}$, we may assume that $a \in N \vee (N' \cap M)$. So, we may assume that $a = xy$ with $x \in N' \cap M$ and $y \in N$. Because $N$ is a factor, $E_{N' \cap M}(z) = \tau(z) 1$ for all $z \in N$. Hence, $E_{N' \cap M}(xyw_n) = \tau(yw_n) x$ and this last sequence converges to $0$ in $\|\cdot\|_2$.

{\bf Step 2.} Let $\xi \in L^2(M)$. The sequence $\|E_{N' \cap M}(\xi w_n)\|_2$ converges to $0$, whenever $w_n$ is a bounded sequence in $N$ that converges weakly to zero. This follows immediately from Step 1.

{\bf Step 3, proof of the lemma.} Define $K$ as the closure of $NbN$ in $L^2(M)$. Since $N \subset M$ is quasi-regular, we may assume that $\dim(K_N) < \infty$. We then find $\xi \in \M_{1,n}(\C) \ot K$ and a, possibly non-unital, $^*$-homomorphism $\psi : N \recht \M_n(\C) \ot N$, such that $x \xi = \xi \psi(x)$ for all $x \in N$ and such that $K$ equals the closure of $\xi (\M_{n,1}(\C) \ot N)$. So, we may assume that $b = \xi d$ for some $d \in
\M_{n,1}(\C) \ot N$. But then, $a w_n b = a \xi \psi(w_n)d$. Since $\psi(w_n) d$ is a bounded sequence in $\M_{n,1}(\C) \ot N$ that converges weakly to zero, the lemma follows from Step 2.
\end{proof}

We are now ready to prove Proposition \ref{prop.unitary}, using an ultrapower argument.

\begin{proof}[Proof of Proposition \ref{prop.unitary}]
Let $N \subset (M,\tau)$ be a quasi-regular inclusion. Suppose that $N$ is a II$_1$ factor.

{\bf Claim 1.} Let $\om$ be a free ultrafilter on $\N$ and $f \in N^\om$ a non-zero projection. If $\cV \subset M^\om$ is a countable set with $E_{(N' \cap M)^\om}(fxf) = 0$ for all $x \in \cV$, there exists a non-zero partial isometry $v \in fN^\om f$ satisfying $vv^* = v^*v$ and $E_{(N' \cap M)^\om}(y) = 0$ for every product $y$ with factors alternatingly from $\cV$ and $\{v^k \mid k \in \Z, k \neq 0\}$.

{\bf Claim 2.} Let $\om$ be a free ultrafilter on $\N$ and $\cV \subset M^\om$ a countable set with $E_{(N' \cap M)^\om}(x) = 0$ for all $x \in \cV$. There exists a unitary $u \in N^\om$ satisfying $E_{(N' \cap M)^\om}(y) = 0$ for every product $y$ with factors alternatingly from $\cV$ and $\{u^k \mid k \in \Z, k \neq 0\}$.

{\bf Proof of Claim 1.} Write $f=(f_n)$ where $f_n$ is a non-zero projection in $N$ for every $n$. Write $\cV = \{ x_k  \mid k \in \N \}$ and choose representatives $x_k = (x_{k,n})_n$ such that $E_{N' \cap M}(f_n x_{k,n} f_n) = 0$ for all $k,n$.  By Lemma \ref{lem.partial}, take partial isometries $v_n \in f_n N f_n$ such that $v_n v_n^* = v_n^* v_n$, $\tau(v_n v_n^*) \geq \tau(f_n) / 4$ and $\|E_{N' \cap M}(y)\|_2 < 1/n$ whenever $y$ is a product of at most $2n+1$ factors alternatingly from $\{x_{0,n},\ldots,x_{n,n}\}$ and $\{v_n^k \mid 1 \leq |k| \leq n \}$. Then, $v:=(v_n)$ does the job.

{\bf Proof of Claim 2.} Define $I$ as the set of partial isometries $v \in N^\om$ satisfying $vv^* = v^* v$ and $E_{(N' \cap M)^\om}(y) = 0$ whenever $y$ is a product with factors alternatingly from $\cV$ and $\{v^k \mid k \in \Z, k \neq 0\}$. By Zorn's lemma, $I$ admits a maximal element $v$. If $v$ is a unitary, we are done. Otherwise, $vv^*=p < 1$ and we set $f = 1-p$. Define $\cW$ as the (countable) set of products $y$ with factors alternatingly from $\cV$ and $\{v^k \mid k \in \Z, k \neq 0 \}$ and such that the product $y$ starts and ends with a factor from $\cV$. Observe that $E_{(N' \cap M)^\om}(fyf) = 0$ for all $y \in \cW$. Indeed,
$$E_{(N' \cap M)^\om}(fyf) = E_{(N' \cap M)^\om}(y) - E_{(N' \cap M)^\om}(yp) = 0 - E_{(N' \cap M)^\om}(v y v^*) = 0 ; .$$
Using Claim 1, take a non-zero partial isometry $w \in f N^\om f$ satisfying $ww^* = w^* w$ and $E_{(N' \cap M)^\om}(y) = 0$ for every product $y$ with factors alternatingly from $\cW$ and $\{w^k \mid k \in \Z, k \neq 0\}$. Then, $v+w \in I$, contradicting the maximality of $v$.

{\bf Proof of the Proposition.} Consider $\cV \subset M \subset M^\om$ with $E_{N' \cap M}(x) = 0$ for all $x \in \cV$. Claim 2 yields a unitary $u \in N^\om$ such that $E_{(N' \cap M)^\om}(y)=0$ for every product $y$ with factors alternatingly from $\cV$ and $\{u^k \mid k \in \Z, k \neq 0\}$. Writing $u = (u_n)$ with $u_n$ unitary for all $n$, some $u_n$ for $n$ big enough will do the job since the elements of $\cV$ are represented by constant sequences in $M^\om$.
\end{proof}

\section{Appendix. Intertwining bimodules and quasi-normalizers} \label{sec.intertwining}

We briefly recall Popa's technique of intertwining subalgebras of a II$_1$ factor using bimodules, introduced in \cite{P1,P5} (see also Appendix C in \cite{V}).

\begin{definition}
Let $(M,\tau)$ be a von Neumann algebra with faithful normal tracial state $\tau$. Suppose that $A,B \subset M$ are von Neumann subalgebras. We say that \emph{$A$ embeds into $B$ inside $M$} and write $A \embed{M} B$, if one of the following equivalent conditions is satisfied.
\begin{itemize}
\item $L^2(M,\tau)$ admits a non-zero $A$-$B$-subbimodule $H$ satisfying $\dim(H_B) < \infty$.
\item $\la M,e_B \ra^+ \cap A'$ contains an element $x$ with $0 < \Tr(x) < \infty$.
\item There exists a projection $p \in B^n$, a normal $^*$-homomorphism $\psi : A \recht p B^n p$ and a non-zero partial isometry $v \in \M_{1,n}(\C) \ot M$ satisfying $x v = v \psi(x)$ for all $x \in A$.
\item There does not exist a generalized sequence $(u_i)_{i \in I}$ of unitaries in $A$ satisfying
$$\|E_B(a u_i b)\|_2 \recht 0 \quad\text{for all}\;\; a,b \in M \; .$$
\end{itemize}
We write $A \fembed{M} B$, if one of the following equivalent conditions is satisfied.
\begin{itemize}
\item For every non-zero projection $p \in M \cap A'$, $L^2(pM,\tau)$ admits a non-zero $Ap$-$B$-subbimodule $H$ satisfying $\dim(H_B) < \infty$.
\item For every $\eps > 0$, there exists a projection $p \in B^n$, a normal $^*$-homomorphism $\psi : A \recht p B^n p$ and a partial isometry $v \in \M_{1,n}(\C) \ot M$ satisfying $\tau(1-vv^*) < \eps$ and $x v = v \psi(x)$ for all $x \in A$.
\end{itemize}
\end{definition}

Let $A \subset (M,\tau)$. In the Preliminaries section, the set $\QN_M(A)$ of elements quasi-normalizing $A$ was introduced, as well as the quasi-normalizer $\QN_M(A)\dpr$. Then, $\QN_M(A)\dpr$ is as well the weak closure of all $x \in M$ for which the closure of $AxA$ in $L^2(M,\tau)$ has finite dimension both as a right $A$-module and as a left $A$-module.

Let $A,B \subset (M,\tau)$. Define
$$p = \vee \{p_0 \mid p_0 \in \la M,e_B \ra \cap A' \;\;\text{is a projection satisfying}\;\; \Tr(p_0) < \infty \} \; .$$
Then, $p L^2(M,\tau)$ clearly is an $A$-$B$-subbimodule of $L^2(M,\tau)$. In fact, it is easy to check that it actually is a $\QN_M(A)\dpr$-$\QN_M(B)\dpr$-subbimodule.

\end{document}